\documentclass[11pt]{amsart}

\newtheorem{theorem}{Theorem}[section]
\newtheorem{lemma}[theorem]{Lemma}

\theoremstyle{remark}

\begin{document}

\title{On punctual Quot schemes for algebraic surfaces.}
\author{V.Baranovsky}
\address{Department of Mathematics, University of Chicago
Chicago, IL 60637}
\email{barashek@math.uchicago.edu}
\date{}

\maketitle

\newcommand{\C}{\mathbb{C}}
\newcommand{\Oo}{\mathcal{O}}
\newcommand{\Or}{\mathcal{O}^{\oplus r}}
\newcommand{\OR}{\mathcal{O}^{\oplus (r-1)}}
\newcommand{\Ox}{\mathcal{O}_{\xi}}
\newcommand{\Q}{\mathrm{Quot}_{[s]}}
\newcommand{\Hi}{\mathrm{Hilb}_{[s]}}
\newcommand{\HP}{\mathrm{Hilb}^d(\mathbb{P}^2)}
\newcommand{\Hom}{\mathrm{Hom}}
\newcommand{\Ext}{\mathrm{Ext}}

\section*{Introduction.}

Let $S$  be a smooth  projective  surface  over  the field of  complex
numbers $\C$. Fix a  closed point $s  \in S$  and  a pair of  positive
integers $r,  d$. By results  of Grothendieck (cf. \cite{G}, \cite{S})
there  exists a projective  scheme $\Q(r,   d)$ parametrizing all quotient
sheaves $\Or_S \to A $ of length $d$ supported at $s$.
We consider this scheme with  its \emph{reduced} scheme structure  and
call it the \emph{punctual Quot scheme}.

Note that $\Q(1,d)$ is nothing but the punctual Hilbert scheme $\Hi^d$
studied by Brian\c con and Iarrobino  in \cite{B}, \cite{I}.  The main
result of  this paper is  the following extension  of their results to
the case $r > 1$:

\medskip

\noindent\textbf{Main Theorem.}  $\Q(r, d)$  is an 
irreducible scheme of dimension $(rd-1)$.

\medskip

 We prove this by exibiting a dense open subset
in $\Q(r, d)$  isomorphic to a rank $(r-1)d$ vector bundle
over $\Q(1, d) = \Hi^d$.

One can show that, for  a quotient $\Or \to A$   as above, the  $d$-th
power of  the maximal ideal  $\mathfrak{m}_{S,  s}$ acts trivially  on
$A$.  Hence the  punctual Quot scheme $\Q(r,  d)$  does not  depend on
$S$  and in  our   proof we   can assume  that   $S =\C^2$. In 
this case a straghtforward generalization of
Nakajima's construction  for Hilbert schemes allows to prove the result.

\medskip

\noindent \textbf{Remark.} The original results of Brian\c con and Iarrobino
were used by G\"ottsche and Soergel in 
\cite{GS} to show that the natural map $\pi: Hilb^d(S) 
\to Sym^d(S)$ is strictly semismall with respect to the natural 
stratifications. This leads to  a simple proof of G\"ottsche's
formula for the Poincar\'e polynomials of $Hilb^d(S)$. Similarly, the
Main Theorem above can be used to show that the natural map $
\pi: M^G(r, d) \to M^U(r, d)$ from the Gieseker moduli space of stable
rank $r$ 
sheaves to the Uhlenbeck compactification of the instanton moduli space, 
is also strictly semismall
(at least in the coprime and unobstructed case). This allows one to find a
connection between  some homological invariants of these spaces. 
A systematic treatment of this questions will appear in the author's
forthcoming paper.

\bigskip

\noindent\textbf{Acknowledgments.} This work was originally motivated
by a conjecture due to V. Ginzburg on semismallness of the map
$\pi: M^G(r, n) \to M^U(r, n)$. The author thanks him for providing
the motivation and also for his helpful discussions and support.
The author also thanks J.Li who explained the role of the
punctual Quot scheme (and also  suggested  an alternative proof of the
irreducibility  statement).

\bigskip

After the first  draft of this  paper was finished, the author learned
about a  preprint by G.Ellingsrud  and M.Lehn \cite{EL} who prove the
same result (among others) by a different method.

\section{Punctual Hilbert scheme.}

The result of this section is well known (cf. \cite{B}, \cite{BI},
\cite{I}). The outline of the proof is given here for convenience of the
reader. It is a slight modification of Corollary 1.2 in \cite{ES}.

\begin{theorem}    $\Hi^d$ is    irreducible  of dimension $(d-1)$.
\end{theorem}

\begin{proof}  First    of  all,     we  can assume      that
  $S=\mathbb{P}^2$.  There  exists  a $\C^*$-action on  $\mathbb{P}^2$
  such    that  our  point $s$   is  a   zero-dimensional  cell of the
  corresponding  Bialynicki-Birula   decomposition.   It  follows that
  $\Hi^d$ is stable  under  the induced  $\C^*$-action on  the  global
  Hilbert scheme $\HP$ of points on  the projective plane. Recall that
  $\HP$   is  \emph{smooth}, (cf.  \cite{N})   and  hence one  has the
  Bialynicki-Birula decomposition  for the  torus action. Then $\Hi^d$
  is  a union  of  cells of  this  decomposition.  One  can prove that
  $\Hi^d$  has a unique $(d-1)$-dimensional  cell and no cells of higher
  dimension (cf. \cite{ES}). Hence $\mathrm{dim}(\Hi^d)= d-1.$
 
  To  prove the   irreducibility  of $\Hi^d$  consider the   universal
  subscheme $Z \subset  \HP\times \mathbb{P}^2$.  By definition of the
  Hilbert   scheme $Z$  is  finite  and  flat over   $\HP$.  Denote by
  $Z_{d-1}$ the subscheme of all points in $Z$ where $d$ sheets of the
  map $f: Z \to \HP$ come together (i.e.   $Z_{d-1}$ is the $(d-1)$-st
  ramification locus of $f$, cf. \cite{GL} for a rigorous definition).
  Then     $(Z_{d-1})_{red}$ is     a  locally    trivial  bundle over
  $\mathbb{P}^2$  with fibers  isomorphic  to $\Hi^d$.   Since $Z$  is
  normal   (cf.   \cite{F}) and  $\HP$ is   smooth,  we  can apply the
  following result due to Lazarsfeld (cf.  \cite{GL} for the statement
  of the   result, the proof of   it is contained  in Lazarsfeld's PhD
  thesis):
 
\emph{Let $f: Z \to H$ be  a finite surjective morphism of irreducible
  varieties, with  $Z$  normal  and  $H$  non-singular. If $Z_{d-1}$ is
  not empty ,  then  every   irreducible   component of   $Z_{d-1}$  has
  codimension $\leq(d-1)$ in $Z$.}

It follows  that any  irreducible component  of  $\Hi^d$ should  be at
least   $(d-1)$-dimensional. Since  $\Hi^d$   has    only one
$(d-1)$-dimensional cell, it can  have only one irreducible component.
\end{proof}

\section{Proof of the Main Theorem.}

Our strategy is   to  find a dense irreducible open subset 
$W \subset \Q(r, d)$ of dimension $(rd-1)$. 

  We define $W$ as the set of all quotients $\Or 
\stackrel{\phi}\longrightarrow A$, $\phi=(\phi_1+ \phi_2 + \ldots + \phi_r)$
such that the first component $\phi_1: \Oo \to A$ is surjective (this
is clearly an open condition). Such a $\phi_1$ corresponds to a point
in $\Hi^d$. Once $\phi_1$ is chosen, the other components
$(\phi_2, \ldots, \phi_r)$ are given by an arbitrary element of
$Hom(\Oo^{\oplus(r-1)}, A) = \C^{(r-1)d}$. Therefore $W$ is a rank
$(r-1)d$ vector bundle over $\Hi^d$. By results of
Brian\c con and Iarrobino, $W$ is irreducible of dimension $(rd-1)$. 

  Now we want to show that $W$ is dense in $\Q(r, d)$. In fact, for
any point $x \in \Q(r, d)$ we will find an irreducible rational curve
$C \subset \Q(r, d)$ connecting it with some point in $W$.

To that end,  we generalize Nakajima's  construction (cf. \cite{N}) of  
the global Hilbert  scheme  $\mathrm{Hilb}^d(\mathbb{C}^2)$ to the 
Quot scheme.   Once we  do that, the   existence  of the   
irreducible curve  will amount  to  an exercise in linear algebra 
(cf. Lemma 2.3).

\bigskip

 Fix a  complex  vector space $V$ of   dimension $d$, and  $N_d$ let 
be the space of pairs   of commuting nilpotent  operators  on $V$.  The space
$N_d$  is naturally a  closed    affine subvariety of   $\mathrm{End}(V)  
\oplus \mathrm{End}(V)$.    

  Consider a subspace $U_r$   of $N_d \times V^{\oplus r}$ formed 
by  all  elements $(B_1,  B_2, v_1, \ldots,  v_r)$ such that there is 
no proper subspace of  $V$ which is invariant under $B_1, B_2$ and 
contains $v_1, \ldots, v_r$. Then $U_1 \times V^{\oplus (r-1)} \subset
U_2 \times V^{\oplus (r-2)} \subset \ldots \subset U_r$ is a chain
of open subsets in $N_d \times V^{\oplus r}$ (each of them is given by 
a condition saying that some system of vectors in $V$ has maximal rank).

\medskip

 Note thate the general linear  group $GL(V)$ acts naturally  
on $V_r$ and it is easy
to prove that $U_r$ is $GL(V)$-stable.

\begin{lemma} 
$GL(V)$ acts freely on $U_r$. 
\end{lemma}

\begin{proof}   Suppose    $g \in  GL(V)$ stabilizes
$(B_1, B_2,   v_1, \ldots, v_r)   \in  U_r$. Then  $\mathrm{Ker}(1-g)$
contains $ v_1, \ldots, v_r$. Since it is also preserved by $B_1, B_2$
, we have $Ker(1-g)=V$ and therefore $g=1$.
\end{proof}

 The following lemma  gives  an explicit construction of  the
punctual Quot scheme:

\begin{lemma} There exists a morphism $\pi: U_r \to \Q(r, d)$ such that

(i) $\pi$ is surjective;

(ii) the fibers of $\pi$ are precisely the orbits of $GL(V)$ action on 
$U_r$;

(iii) $\pi^{-1}(W) = U_1 \times V^{\oplus (r-1)}$. 
\end{lemma}

\begin{proof} We can  assume that $S= \C^2  = Spec \;\C[x_1, x_2]$ and $s=0
  \in \C^2$.

To  construct $\pi$ suppose   that $(B_1, B_2, v_1, \ldots,
v_r)$  is a  point  in  $U_r$ and   consider  a $\C[x_1,  x_2]$-module
structure on $V$ in which $x_1$ acts by  $B_1$ and $x_2$ acts by $B_2$. 
We can view $V$ as a quotient of a free $\C[x_1,   x_2]$-module with 
generators  $v_1, \ldots,  v_r$.  Since $B_1$ and $B_2$ are  nilpotent 
$\sqrt{\mathrm{Ann}(V)} = (x_1, x_2)$. Therefore a  coherent sheaf $A$  
on  $\C^2$ associated  with $V$ is  a quotient of $\Or$ supported 
at $s$. Moreover,  $V \simeq \mathrm{H}^0(S, A)$ as vector spaces.

A different point in the same $GL(V)$-orbit  defines 
 an isomorphic  quotient, hence the fibers of $\pi:  U_r/ \mathrm{GL(V)} 
\to \Q(r,d)$  are $GL(V)$-invariant. Moreover, suppose that
 two points $u_1, u_2$ 
of  $U_r$ give rise to isomorphic quotients $A_1$, $A_2$. Then the 
induced isomorphism between $\mathrm{H}^0(S, A_1)$ and $\mathrm{H}^0(S, A_2)$
defines an element of $GL(V)$ taking $u_1$ to $u_2$. Therefore,
each fiber of $\pi$ is precisely one $GL(V)$-orbit. This proves $(ii)$.

To prove $(i)$, suppose we have a quotient $ \Or \to  A \to 0 $ of length
  $d$ supported  at zero. Multiplication  by $x_1$ and $x_2$ induces a
  pair of commuting nilpotent operators on $\mathrm{H}^0(S, A)$. Choose a
  $\C$-linear isomorphism  $\mathrm{H}^0(S, A) \simeq  V$. The generators
  of the free  $\C[x_1,  x_2]$-module $\mathrm{H}^0(S, \Or)$  project  to
  some vectors $v_1,  \ldots, v_r$ in  $V$.  Since $v_1,  \ldots, v_r$
  generate  $V$  as a $\C[x_1, x_2]$-module,   $(x_1, x_2, v_1, \ldots
  v_r)$ is a point of $U_r$. Thus $(i)$ is proved.

  Finally, $(iii)$ follows from definitions of $W$ and $U_1$.
\end{proof}
  
  Now we want to show that any point in $U_r$ can be deformed
to a point in the preimage of $W$. The above construction will allow us to
construct this deformation using the following lemma

\begin{lemma}
Let $B_1$, $B_2$ be two commuting nilpotent operators on a vector space V.
There exists a third nilpotent operator $B_2'$ and a vector $w \in V$
such that

(i) $B_2'$ commutes with $B_1$;

(ii) any linear combination $\alpha B_2 + \beta B_2'$ is nilpotent;

(iii) $(B_1, B_2', w) \in U_1$, \;i.e. $w$  is a cyclic vector for the
pair of operators $(B_1, B_2')$.
\end{lemma}

  This lemma will be proved later. Now we will show how it can be used
to give a

\bigskip 

\noindent\emph{Proof of the Main Theorem:} 

Let $x$ be a point of $\Q(r, d)$ and $u_1 = (B_1, B_2, v_1, \ldots, v_r)$
be any point of $\pi^{-1}(x) \subset   U_r$. Choose a nilpotent
operator $B_2'$ and a vector $w \in V$ as in Lemma 2.3. 
Connect the points $u_1$ and $u_2 = (B_1, B_2', w, v_2 , \ldots, v_r)$ with
a straight line $\Phi(t)$, $t \in \C$ \; such that  $\Phi(1) = u_1$ and 
$\Phi(0) = u_2$. This $\Phi(t)$ is given by  equation:
$$
\Phi(t)= (B_1, tB_2' + (1-t)B_2, tw +(1-t)v_1, v_2, \ldots, v_r)
$$

  Note  that  for all $t \in C$,\;   $B_2(t) = tB_2' + (1-t)B_2 $   
is nilpotent and commutes with $B_1$. Therefore the image of $\Phi(t)$
is a subset of $N_d \times V^{\oplus r}$. Since $U_r$ is open in 
$N_d \times V^{\oplus r}$, there is a dense open subset 
$C \subset \C$ such that $\Phi(C) \subset U_r$.
Similarly, there exists a dense open subset $C_1 \subset 
C$ such that $\Phi(C_1) \subset U_1 \times V^{\oplus 
(r-1)}$.

  Hence the image $\pi(\Phi(C)) \subset \Q(r, d)$ is an 
irreducible rational curve connecting $x = \pi(u_1)$ with $\pi(u_2) \in W$.
Note that $\pi(\Phi(C_1)) \subset W$.
Therefore $x$ belongs to the closure of $W$. Since by Theorem 1.1
$W$ is irreducible of dimension $(rd-1)$, the scheme $\Q(r, d)$ is
also irreducible of dimension $(rd-1)$. The Main Theorem is proved.

\bigskip 

\noindent\emph{Proof of Lemma 2.3:}

\medskip 

\emph{Step 1.}
 We will find a basis $e_{i,j}$ of $V$, where $1 \leq i \leq k$ and
$1 \leq j \leq \mu_i$ such that

(a) $B_1^{j-1}(e_{i, 1}) = e_{i, j}$ for $j \leq \mu_i$ and 
$B_1^{\mu_j}(e_{i, 1}) = 0$ (i.e. $B_1$ has Jordan canonical form
in the basis $e_{i, j}$);

(b) $B_2(e_{i, 1}) \in \Big(\bigoplus_{k \geq i +1} \C\; e_{k, 1}) \oplus
B_1 \cdot V$.

\bigskip 

 To that end, recall one way to construct a Jordan basis for $B_1$.
Let $d = \dim V$ and $V_i = Ker(B_1^{d-i})$. The subspaces $V_i$ form 
a decreasing filtration $V = V_0 \supset V_1 \supset V_2 \ldots $.
Moreover, $B_1 \cdot V_i \subset V_{i+1}$. Firstly, we choose a basis
$(w_1, \ldots, w_{a_1})$ of $W_1: = V_0/V_1$. Lift this basis to some 
vectors $e_{1, 1}, e_{2, 1}, \ldots, e_{a_1, 1}$ in $V_0$ and set
all $\mu_1, \ldots, \mu_{a_1}$ equal to $d$. Secondly, choose a basis 
$(w_{a_1 + 1}, \ldots, w_{a_2})$ of $W_2: = V_1/ (B_1\cdot V_0 + V_2)$. 
Lift this basis to some  vectors $e_{a_1 + 1, 1}, e_{2, 1}, \ldots, 
e_{a_2, 1}$ in $V_1$ and set all $\mu_{a_1+1}, \ldots, \mu_{a_2}$ equal 
to $d-1$. Continue in this manner by choosing bases of the spaces
$W_{i+1} = V_i/ (B_1 \cdot V_{i-1} + V_{i+1})$ and lifting them to $V_i$.
This procedure gives us  vectors $e_{1, 1}, e_{2, 1}, \ldots,
e_{k, 1}$ and the formula (a) tells us how to define $e_{i, j}$ for
$j \geq 2$. It is easy to check that the system of vectors $\{e_{i, j}\}$
is in fact a basis of $V$.

\bigskip 

  If we want to have the property (b)  we should be more careful with 
the choice of $w_i$.  Note that all the subspaces $V_i$ and $B_1 \cdot V_i$
are $B_2$-invariant. Therefore we have  an induced action of $B_2$ on
each $W_i$. We can choose our basis $(w_{a_{i-1} + 1}, \ldots, w_{a_i})$ 
of $W_i$ in such a way that $B_2(w_i) \in \bigoplus_{s=i+1}^{a_i} \C \; w_s$
for all $i \in \{a_{i-1} + 1, \ldots, a_i \}$.
This ensures that $(b)$ holds as well.

\bigskip

\emph{Step 2.}
Define $B_2'$ by $B_2'(e_{i, j}) = e_{i+1, j}$ if 
$j \leq \mu_{i+1}$ and 0 otherwise. It is immediate that $B_2'$ is nilpotent
and that $[B_1, B_2']=0$. Let $w = e_{1, 1}$. Then $e_{i, j} = B_1^{j-1}
(B_2^{i-1}(w))$ hence $(B_1, B_2', w) \in U_1$.

\bigskip

\emph{Step 3.}
Note that both $B_2$ and $B_2$ are lower-triangular with zeros on the 
diagonal in the basis of $V$ given by
$$
e_{1, 1}, e_{2, 1}, \ldots, e_{k, 1}, e_{1, 2}, e_{2, 2}, \ldots
 e_{k, 2}, \ldots
$$ 
Hence any linear combination of $B_2$ and $B_2'$ as also lower-triangular
and has zeros on the diagonal. Therefore $\alpha B_2 + \beta B_2'$ is
nilpotent for any complex $\alpha$ and $\beta$. This completes the proof
of Lemma 2.3.

\bibliographystyle{amsplain}

\end{document}